\documentclass[letterpaper, 10 pt, conference]{ieeeconf}
\IEEEoverridecommandlockouts

\usepackage{amsmath}
\usepackage{graphicx}
\usepackage{cite}
\usepackage{authblk}
\usepackage{amssymb}
\usepackage{subfigure}
\usepackage{multirow}
\usepackage{color}
\usepackage{tikz} \usetikzlibrary{arrows,calc,decorations.pathreplacing}

\newtheorem{theorem}{Theorem}

\newtheorem{definition}{Definition}
\newtheorem{corollary}{Corollary}

\title{Optimal Strategies for the Game of Protecting a Plane in 3-D}
\author{Eloy Garcia, Isaac Weintraub, David W. Casbeer, and Meir Pachter
\thanks{This work has been supported in part by AFOSR LRIR No. 21RQCOR084.}
\thanks{E. Garcia, I. Weintraub, and D. Casbeer are with the Control Science Center of Excellence, Air Force Research Laboratory, Wright-Patterson AFB, OH 45433. Corresponding author \ttfamily{eloy.garcia.2@us.af.mil}}
\thanks{M. Pachter is with the Department of Electrical Engineering, Air Force Institute of Technology, Wright-Patterson AFB, OH 45433.}
}

\begin{document}
\maketitle 
\begin{abstract}
A conflict between rational and autonomous agents is considered. The paper addresses a differential game of protecting a target in the 3-D space. This problem highlights the strong correlation between the highly dynamic scenario, the uncertainty on the behavior of the adversary, and the online and robust computation of state-feedback strategies which guarantee the required level of performance of each player. 
This work significantly expands previous results around this problem by providing the players' state-feedback saddle-point strategies. Additionally, the continuously differentiable Value function of the multi-agent differential game is obtained and it is shown to be the solution of the Hamilton-Jacobi-Isaacs equation. 
Finally, the Barrier surface is explicitly obtained and illustrative examples highlight the robustness properties and the guarantees provided by the saddle-point strategies obtained in this paper.
\end{abstract}

\section{Introduction} \label{sec:intro}
Multi-agent dynamic conflicts are a representative problem showcasing autonomous systems, which are deployed in uncertain environments. These problems are generally addressed by means of differential game theory, which provides a framework to analyze pursuit-evasion scenarios, reach-avoid games, and conflicts between adversarial teams. The foundations of differential game theory can be found in the seminal work by Isaacs \cite{Isaacs65}. In general, pursuit-evasion games consider an evader trying to escape from a pursuer without intentionally aiming at reaching a particular region of the game set  \cite{breakwell1979point,Bopardikar09,liu2014evasion,takei2014efficient}.
More recently, reach-avoid games \cite{chen2014path,zhou2018efficient,Blake04} and attack-defend problems \cite{Coon17,Garcia20Multi,fisac2015pursuit,Fuchs17,Garcia2019}  address more general conflicts where the evader aims at reaching a target or goal set while also avoiding the pursuer. Therefore, two distinct outcomes exist, the evader is able to reach the protected target or the evader is intercepted by the pursuer before reaching said target. 

 The differential game framework is generally desired in order to solve pursuit-evasion, reach-avoid, and attack-defend games  \cite{chen2016multiplayer,weintraub20intro}, but it is often avoided due to the challenges in solving the Hamilton-Jacobi-Isaacs (HJI) equation \cite{pierson2016intercepting}.
This paper follows Isaacs' approach and enables cooperation between players of the same team. 
Most of the references addressing pursuit-evasion and reach-avoid games assume that the agents play in the 2-D Cartesian plane. 
 This paper extends the differential game approach to address reach-avoid scenarios in the 3-D Cartesian space.  Additionally, the differential game is solved in analytical and closed-form. The obtained strategies  can be implemented on-line in order to take advantage of non-optimal behaviors by the opponent.


Several papers have considered reach-avoid games in the  3-D Cartesian space. The paper \cite{Garcia20RAL} addressed an extension of the capture-the-flag game in 3-D with one evader and two cooperative pursuers.
Reference \cite{Weintraub2020} considered the defense of a non-maneuverable target in the 3-D Cartesian space with one pursuer and one evader.  The pursuer was assumed to be faster than the defender. In the present paper, the pursuers and the evader have equal speed and it is shown that, in the case of equal speed players, at least three cooperative pursuers are needed to intercept the evader. A shortcoming of \cite{Weintraub2020} is that proposed strategies  were not verified. Obtaining the solution via Isaacs' method is the ideal situation in differential games, as it provides guarantees of correctness of the solution  \cite{chen2016multiplayer}. Verification is often overlooked due to challenges in solving the HJI equation; however, verification is of great importance in order to guarantee that the unique saddle-point solution has been in fact synthesized. In addition, we provide the full cooperative strategy. This aspect is usually overlooked in many multi-player games since approximation of the solutions are based on decomposition of strategies which severely limit the level of cooperation of a given team.

Recently, the authors of  \cite{yan2019construction} considered the same problem as in this paper but they only addressed the Game of Kind by determining the Barrier surface. The Game of Kind (which player wins) and the Game of Degree (how the player wins) are equally important and intrinsically related. The Game of Kind solution is necessary for the players to determine which Game of Degree to play and the solution to the Game of Degree needs to sustain the Barrier surface obtained by solving the Game of Kind. The latter means, in technical terms, that optimal strategies should result in a semipermeable Barrier surface. 

In this paper we solve the Game of Kind in a more explicit form than   \cite{yan2019construction} and, more importantly, we synthesize and verify the saddle-point state-feedback strategies, thus, providing a complete solution to the Game of Degree as well. In other words, finding the winning regions based on initial conditions is not useful by itself if the players are not able to actually realize such outcomes.
Another main difference with respect to \cite{yan2019construction} is that we consider state-feedback strategies. State-feedback strategies represent a more general and practical class of strategies than the class of strategies considered in \cite{yan2019construction}, where the evader needs to communicate its current strategy to the adversary, the group of pursuers. State-feedback strategies, on the other hand, are only based on the current state of the system and the players do not need to know the current strategy implemented by the adversary. This is more practical and easier to implement in actual conflicts where, for obvious reasons, the players will not share strategic information with the opponent. One can imagine, for instance, that the evader will deceive its adversary by sharing false strategic information.

When considering state-feedback strategies, we have that at least three pursuers are needed to capture the evader. By considering state-feedback strategies, we are able to solve both the Game of Degree and the Game of Kind. These solutions represent a valuable contribution with respect to \cite{yan2019construction} where only the Game of Kind was addressed. The solution of the Game of Kind only answers the question of which team wins but it does not provide information on how the players can achieve such outcome. 
It is very important to solve the Game of Degree in order to provide the complete solution of the differential game. By solving the Game of Degree, we are able to determine the optimal strategies of all players and actually enforce the outcome predicated by the Game of Kind solution. In this paper, we obtain and verify this solution. The Value function is obtained and it is shown to be continuously differentiable and to be the solution of the HJI equation.

The rest of the paper is organized as follows. The problem is formulated within the differential game framework in Section \ref{sec:Prelim}. The state-feedback, optimal strategies are derived and verified in Section \ref{sec:1P1E}. The Barrier surface of this reach-avoid game in 3-D is explicitly obtained in Section \ref{sec:Barrier}.
Illustrative examples are shown in Section \ref{sec:ex} and conclusions are drawn in Section \ref{sec:concl}.


\section{The Differential Game in 3D}   \label{sec:Prelim}
Consider a reach-avoid game with three pursuers and one evader, where all players have the same speed. The pursuers are denoted by $P_1$, $P_2$, and $P_3$. The group of pursuers cooperate in order to capture the evader, denoted by $E$, while defending the goal plane, denoted by $\Omega_G$. The evader aims at reaching the goal plane and,  without loss of generality, we assume that the goal plane, $\Omega_G$, is given by $z(x,y)=0$. The game is played in the Euclidean 3-D space. In particular, the game set is $\Omega=\{x,y,z \ | \ z>0 \}$.

The states of the pursuers are given by their Cartesian coordinates $\textbf{x}_i=(x_i,y_i,z_i)$, for $i=1,2,3$. The state of $E$ is specified by its Cartesian coordinates $\textbf{x}_E=(x_E,y_E,z_E)$. Similarly, 
The complete state of the differential game is defined by $\textbf{x}:=( \textbf{x}_E, \textbf{x}_1, \textbf{x}_2,\textbf{x}_3)\in \mathbb{R}^{12}$. The control set, for each agent, is defined as  $\mathcal{U}:=\{ \textbf{u}\in \mathbb{R}^3 | \Vert \textbf{u} \Vert _2 =1 \}$. The control input of $E$ is denoted by $\textbf{u}_E=(u_{x_E},u_{y_E},u_{z_E})$. The control input of pursuer $i$ is denoted by $\textbf{u}_i=(u_{x_i},u_{y_i},u_{z_i})$, where $\textbf{u}_E,\textbf{u}_i \in \mathcal{U}$, for $i=1,2,3$. Hence, the dynamics/kinematics $\dot{\textbf{x}}=\textbf{f}(\textbf{x},\textbf{u}_E,\textbf{u}_i)$ are specified by the system of ordinary differential equations 
\begin{equation}
\begin{alignedat}{2}
	\dot{x}_E&= u_{x_E},  \qquad  \dot{y}_E=u_{y_E},  \qquad  \dot{z}_E=u_{z_E}  \\
	\dot{x}_i&= u_{x_i},  \qquad  \  \dot{y}_i=u_{y_i},  \qquad \ \ \dot{z}_i=u_{z_i}    \label{eq:xD}
\end{alignedat}
\end{equation}
with $x_E(0)=x_{E_0}, y_E(0)=y_{E_0}, z_E(0)=z_{E_0}, x_i(0)=x_{i_0}, y_i(0)=y_{i_0}, z_i(0)=z_{i_0}$, for $i=1,2,3$. Without loss of generality, the speeds have been normalized. 
The initial state of the system is $\textbf{x}_0 := (\textbf{x}_{E_0},\textbf{x}_{1_0},\textbf{x}_{2_0},\textbf{x}_{3_0}) = \textbf{x}(0)$. It is assumed that $z_{E_0}>0$, that is, the evader is initially located in the game set $\Omega$. The single integrator dynamics are typical in many games by Isaacs \cite{Isaacs65}. It is possible to obtain explicit and closed-loop strategies which can be implemented online. This is of great importance since the agents are able to cooperate at the highest level and to take advantage of non-optimal plays by the opponent. In addition, these strategies can be used as a robust approximation when the players display more complex dynamics with dynamic constraints such as turning rate constraints and acceleration constraints \cite{GarciaACC21BVR}.

The game under consideration is a two-termination set differential game \cite{getz1979qualitative,Getz1981capturability,ardema1985combat}. One terminal condition is capture of the evader by any of the pursuers. In such a case, the pursuer team wins the game. Alternatively, the evader wins if it can reach the goal plane $\Omega_G$ before being captured. Hence, the termination set is
\begin{align}
   \mathcal{T} :=  \mathcal{T}_p    \   \bigcup \   \mathcal{T}_e   \label{eq:TwoSets}
\end{align}
where 
\begin{align}
\left.
	 \begin{array}{l l}
 \mathcal{T}_p:= \!  \big\{  \textbf{x} \ | \ \| \textbf{x}_1\!-\!\textbf{x}_E  \| _2 =0 \big\} \cup \big\{  \textbf{x} \ | \ \| \textbf{x}_2\!-\!\textbf{x}_E  \| _2 =0 \big\} \\
\qquad \quad  \cup \big\{  \textbf{x} \ | \ \| \textbf{x}_3\!-\!\textbf{x}_E  \| _2 =0 \big\}
  \end{array}   \right.    \label{eq:Termc}
\end{align}
represents the outcome where $E$ is captured before reaching the target and
 \begin{align}
  \mathcal{T}_e:= \big\{  \textbf{x} \ | \  z_E  =0 \big\}     \label{eq:Terme}
\end{align}
represents the outcome where  $E$ wins the game by reaching $\Omega_G$.
The terminal time $t_f$ is the time instant when the state of the system satisfies \eqref{eq:TwoSets}, at which time the terminal state is $\textbf{x}_f: = (\textbf{x}_{E_f},\textbf{x}_{1_f},\textbf{x}_{2_f},\textbf{x}_{3_f})   = \textbf{x}(t_f)$.

The concepts of Game of Kind and Game of Degree are fundamental in differential game theory \cite{Isaacs65}. The solution to the Game of Kind determines which team wins the game. Solving the Game of Degree provides the value of the game and the saddle-point strategies that realize the outcome prescribed by the Game of Kind. 
Because of the two different outcomes specified in  \eqref{eq:TwoSets}, the Game of Kind needs to be solved in order to partition the state space into two winning regions, one for evader and one for the pursuer team. Since different Games of Degree are played in each region, it is essential for each player to determine which region the current state of the system is in. Armed with this information, each player can then implement the appropriate optimal strategy for the corresponding Game of Degree.
 The state space $\mathbb{R}^{12}$ is partitioned into two sets: $\mathcal{R}_p$ and $\mathcal{R}_e$ which are defined as follows
\begin{align}
\left.
	 \begin{array}{l l}
\mathcal{R}_p:=  \big\{ \ \textbf{x} \ |  \ B( \textbf{x})>0  \big\}, \ \
\mathcal{R}_e:=  \big\{ \ \textbf{x} \ |  \  B( \textbf{x})<0  \big\}. 
\end{array}   \right.  \label{eq:ReDef}
\end{align}
The Barrier surface, which separates the two sets $\mathcal{R}_p$ and $\mathcal{R}_e$, is specified by
\begin{align}
\left.
	 \begin{array}{l l}
\mathcal{B}:=  \big\{ \ \textbf{x} \ | \ B( \textbf{x})=0  \big\} 
\end{array}   \right.  \label {eq:BarrSurface}
\end{align}
where the Barrier function, $B( \textbf{x})$, is explicitly obtained in Section \ref{sec:Barrier}. In this paper we consider the Game of Degree when  $\textbf{x} \in \mathcal{R}_p$. The terminal performance functional is
\begin{align}
  J(\textbf{u}_E(t),\textbf{u}_i(t);\textbf{x}_0)=\Phi_p(\textbf{x}(t_f))	\label{eq:costDG}
\end{align}
where  $\Phi_p(\textbf{x}(t_f)):=   z_{E_f}=z_E(t_f)$. The Value of the game is 
\begin{align}
  V(\textbf{x}_0):= \min_{\textbf{u}_E(\cdot)} \ \max_{\textbf{u}_i(\cdot)} J(\textbf{u}_E(\cdot),\textbf{u}_i(\cdot);\textbf{x}_0)	\label{eq:Vsf}
\end{align}
subject to \eqref{eq:xD} and \eqref{eq:Termc}, where $\textbf{u}_E(\cdot)$ and $\textbf{u}_i(\cdot)$ are the teams' state-feedback strategies, for $i=1,2,3$. In other words, the players only have access to the current state of the system; they do not have access to present or future controls of the opponent. The strategies derived in this paper are only a function of the current state of the game.

When the solution of the Game of Kind prescribes that $E$ is going to be captured before reaching the target, $E$ strives to minimize its terminal separation with respect to the goal plane $\Omega_G$ at the time instant of capture. The pursuers aim at intercepting $E$ while maximizing the terminal separation. This strategy provides a practical outcome in case the pursuers do not play optimally. By trying to get as close to $\Omega_G$  as possible, $E$ is in position to exploit a mistake by the pursuers and potentially win the game.

%

Finally, we define the order of preference. Let $\zeta_a$ and $\zeta_b$ be numerical outcomes of terminating plays, where $\zeta_a<\zeta_b$. 

\begin{definition}   \label{def:Order}  
Define $\Upsilon$ as the outcome of any non-terminating play. Then, the differential game when $\textbf{x} \in \mathcal{R}_p$ has an order of preference of Type $E$.
\end{definition}

In a Type $E$ order of preference \cite{Lewin94} the following holds:
\begin{align}
  \zeta_a \stackrel{E}{\succ} \zeta_b \stackrel{E}{\succ} \Upsilon  \ \ \  \text{and} \ \ \  \Upsilon \stackrel{P}{\succ} \zeta_b \stackrel{P}{\succ} \zeta_a. \nonumber
\end{align}
This means that player $E$, who aims at minimizing the numerical outcome of the game, considers non-termination to be inferior to all other outcomes. Player (or team) $P$, who aims at maximizing the numerical outcome of the game, considers non-termination to be superior to all other outcomes. Since the players have the same speed, the evader could avoid capture by  running away from the pursuers and the goal plane. However, from Definition \ref{def:Order}, this strategy will result in the worst outcome for the evader and the best outcome for the pursuer team. Hence, the evader is not advised to follow such that strategy.

\begin{theorem}  \label{th:main}
Consider the differential game \eqref{eq:xD} and \eqref{eq:costDG}. The optimal headings of $E$, $P_1$, $P_2$, and $P_3$ are constant under optimal play and their trajectories are straight lines. 
\end{theorem}
\textit{Proof}. Consider \eqref{eq:xD} and \eqref{eq:costDG},
the optimal control inputs (in terms of the co-state variables) can be immediately obtained from  $\min_{\textbf{u}_E} \max_{\textbf{u}_i} \mathcal{H}$, where the Hamiltonian is
\begin{align}
  \left.
	\begin{array}{l l}
	\mathcal{H}\!\!\! &=  \lambda_{x_E} u_{x_E} +\lambda_{y_E} u_{y_E} +\lambda_{z_E} u_{z_E} \\
	  &~~ + \sum_{i=1}^3 ( \lambda_{x_i} u_{x_i} + \lambda_{y_i} u_{y_i}  +\lambda_{z_i} u_{z_i} ) 
\end{array}  \right.   \nonumber
\end{align}
 where $\lambda^T=( \lambda_{x_E}, \lambda_{y_E}, \lambda_{z_E}, \lambda_{x_i}, \lambda_{y_i}, \lambda_{z_i}) \in \mathbb{R}^{12}$ is the co-state, for $i=1,2,3$. 
Additionally, the co-state dynamics are $\dot{\lambda} = - \frac{\partial \mathcal{H}}{\partial \textbf{x}} = \textbf{0}$. Thus, 
 all co-states are constant and the optimal headings are constant as well. Consequently, the optimal trajectories are straight lines. \ \ \ \ \ \ \ \ \ \ \ \ \ \ \ \ \ \ \ \ \ \ \ \ \ \ \ \ \ \ \ \ \ \  \ \ \ $\square$

\begin{figure}
	\begin{center}
		\includegraphics[width=8.2cm,trim=1.2cm 1.7cm 1.2cm 1.9cm]{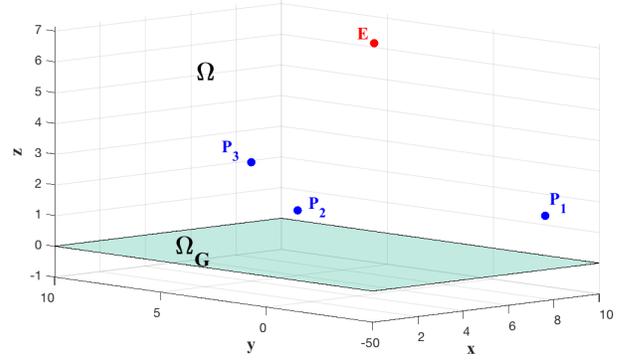}
	\caption{Game of protecting a plane between an evader and a team of three pursuers}
	\label{fig:Problem}
	\end{center}
\end{figure}

\section{Optimal Strategies}  \label{sec:1P1E}

By a change of coordinate we can transform the problem of protecting any plane in the 3-D space to the problem of protecting the plane $z(x,y)=0$ as it is illustrated in Fig. \ref{fig:Problem}.
Let us consider the Game of Degree in the pursuers' winning region. We assume through this section that $\textbf{x} \in \mathcal{R}_p$, equivalently, we assume that $B(\textbf{x})>0$, where $B(\textbf{x})$ is the Barrier function. The Game of Kind is solved in the following section, where the Barrier surface is delineated, that is, the Barrier function, $B(\textbf{x})$, is obtained in closed-form. The Barrier surface separates the winning regions of the evader and the pursuer team.

We define $R_E=\sqrt{x_E^2+y_E^2+z_E^2}$. Also, define $R_i=\sqrt{x_i^2+y_i^2+z_i^2}$, for $i=1,2,3$.
We can characterize the dominance or reachable region of the evader with respect to each pursuer. The reachable regions of $E$ and pursuer $P_i$ are separated by the plane $H_i=0$, which is orthogonal to the segment $\overline{EP_i}$, where
\begin{align}
 \left.
	 \begin{array}{l l}
	H_i \!=\! (x_E\!-\!x_i)x\!+\!(y_E\!-\!y_i)y\!+\!(z_E\!-\!z_i)z-\frac{R_E^2 -R_i^2}{2}  
	\end{array}  \right.  \label{eq:OP-Ei}
\end{align}
for $i=1,2,3$.
Let us define 
\begin{align}
 \left.
	 \begin{array}{l l}
	\mathcal{E}_i=  \{ x,y,z \ | \ H_i  > 0 \}  
   \label{eq:halfplane}
\end{array}  \right.  
\end{align}
for $i=1,2,3$. Then, $E$'s reachable region with respect to all pursuers is given by
\begin{align}
 \left.
	 \begin{array}{l l}
	\mathcal{E}= \cap_i \ \mathcal{E}_i.  
\end{array}  \right.   \label{eq:SpyS}
\end{align}
The following theorem provides the optimal strategies for the case where $\textbf{x}\in\mathcal{R}_p$.


\begin{theorem}   \label{th:1P1E}
Consider the differential game of protecting a plane in the 3-D space and
assume that $\textbf{x}\in\mathcal{R}_p$. The Value function is $C^1$ and it is the solution of the Hamilton-Jacobi-Isaacs (HJI) partial differential equation. The Value function is given by
\begin{align}
 \left.
	 \begin{array}{l l}
	V(\textbf{x}) = \frac{ \nu_E R_E^2 - \nu_1 R_1^2  -\nu_2 R_2^2  -\nu_3 R_3^2 }{2\Lambda}
\end{array}  \right.  \label{eq:ValueFn}
\end{align}
where 
\begin{align}
 \left.
	 \begin{array}{l l}
	\nu_E =  x_1(y_2\!-\!y_3)+x_2(y_3\!-\!y_1)+x_3(y_1\!-\!y_2)  \\
	  \nu_1=x_2(y_3\!-\!y_E)+x_3(y_E\!-\!y_2)+x_E(y_2\!-\!y_3)  \\ 
	 \nu_2=x_3(y_1\!-\!y_E)+x_1(y_E\!-\!y_3)+x_E(y_3\!-\!y_1)  \\
	 \nu_3=x_1(y_2\!-\!y_E)+x_2(y_E\!-\!y_1)+x_E(y_1\!-\!y_2)
\end{array}  \right.  \label{eq:zoterms}
\end{align}
and
\begin{align}
 \left.
	 \begin{array}{l l}
	\Lambda \!\!\!\! &=   x_E(y_3z_1\!-\!y_1z_3) \!+\! y_E(x_1z_3\!-\!x_3z_1)  
	   \!+\!z_E(x_3y_1\!-\!x_1y_3) \\
	 & ~~ - x_2[y_3(z_1\!-\!z_E)+y_1(z_E\!-\!z_3)+y_E(z_3\!-\!z_1)] \\  
	 & ~~ - y_2[x_3(z_E\!-\!z_1)+x_1(z_3\!-\!z_E)+x_E(z_1\!-\!z_3)] \\ 
	  & ~~ - z_2[x_3(y_1\!-\!y_E)+x_1(y_E\!-\!y_3)+x_E(y_3\!-\!y_1)]. 
\end{array}  \right.  \label{eq:Denominator}
\end{align}
The optimal  state-feedback  strategy of the evader is given by 
\begin{align}
 \left.
	 \begin{array}{l l}
	 \textbf{u}^*_E=\frac{1}{d_E}[x^*(\textbf{x})-x_E, \ y^*(\textbf{x})-y_E, \ z^*(\textbf{x})-z_E ]  \\   
\end{array}  \right.    \label{eq:OptimalInputE}
\end{align}
where $d_E=\sqrt{(x^*-x_E)^2 + (y^*-y_E)^2+(z^*-z_E)^2}$.
The optimal  state-feedback  strategies of the pursuers are given by 
\begin{align}
 \left.
	 \begin{array}{l l}
	 \textbf{u}^*_i=\frac{1}{d_i}[x^*(\textbf{x})-x_i, \ y^*(\textbf{x})-y_i, \ z^*(\textbf{x})-z_i ]  \\   
\end{array}  \right.    \label{eq:OptimalInputsP}
\end{align}
where  $d_i=\sqrt{(x^*-x_i)^2 + (y^*-y_i)^2+(z^*-z_i)^2}$, for $i=1,2,3$.
The coordinates of the optimal interception point are
\begin{align}
 \left.
	 \begin{array}{l l}
	x^*(\textbf{x}) \!\!\!\!&= \frac{\mu_E R_E^2 - \mu_1 R_1^2  -\mu_2 R_2^2  -\mu_3 R_3^2}{2\Lambda}
\end{array}  \right.  \label{eq:Xopt}
\end{align}
\begin{align}
 \left.
	 \begin{array}{l l}
	y^*(\textbf{x}) \!\!\!\!&= \frac{\eta_E R_E^2 - \eta_1 R_1^2  -\eta_2 R_2^2  -\eta_3 R_3^2 }{2\Lambda}
\end{array}  \right.  \label{eq:Yopt}
\end{align}
and $z^*(\textbf{x})=V(\textbf{x})$, where
\begin{align}
 \left.
	 \begin{array}{l l}
	\mu_E = y_1(z_2\!-\!z_3)+y_2(z_3\!-\!z_1)+y_3(z_1\!-\!z_2)  \\
	 \mu_1=y_2(z_3\!-\!z_E)+y_3(z_E\!-\!z_2)+y_E(z_2\!-\!z_3)  \\ 
	 \mu_2 = y_3(z_1\!-\!z_E)+y_1(z_E\!-\!z_3)+y_E(z_3\!-\!z_1)  \\
	 \mu_3 = y_1(z_2\!-\!z_E)+y_2(z_E\!-\!z_1)+y_E(z_1\!-\!z_2)
\end{array}  \right.  \label{eq:xoterms}
\end{align}
and
\begin{align}
 \left.
	 \begin{array}{l l}
	\eta_E = x_1(z_3\!-\!z_2)+x_2(z_1\!-\!z_3)+x_3(z_2\!-\!z_1)  \\
	 \eta_1 =x_2(z_E\!-\!z_3)+x_3(z_2\!-\!z_E)+x_E(z_3\!-\!z_2)  \\ 
	 \eta_2 =x_3(z_E\!-\!z_1)+x_1(z_3\!-\!z_E)+x_E(z_1\!-\!z_3)  \\
	\eta_3 = x_1(z_E\!-\!z_2)+x_2(z_1\!-\!z_E)+x_E(z_2\!-\!z_1).
\end{array}  \right.  \label{eq:yoterms}
\end{align}
\end{theorem}
\textit{Proof}. The evader, wishing to minimize its terminal separation with respect to the goal plane $z(x,y)=0$, aims at the point on its reachable region which is the closest to the goal plane. Since $\textbf{x}\in\mathcal{R}_p$, the three planes $H_i$, for $i=1,2,3$, intersect at only one point. In addition, the intersection point given by $I^*=(x^*,y^*,z^*)$ is such that $z^*>0$. Therefore, the coordinate of the intersection point can be obtained by solving the system of three linear equations $H_i=0$ in $(x,y,z)$, for $i=1,2,3$. The solution, written explicitly in terms of the state of the system, $\textbf{x}$, is given by \eqref{eq:Xopt}-\eqref{eq:Yopt}, and \eqref{eq:ValueFn}. The Value function is given by the terminal distance between the evader and the goal plane $z(x,y)=0$. Hence, the Value function is given by $V(\textbf{x})=z^*(\textbf{x})$ as it is shown in \eqref{eq:ValueFn}.

We can now compute the gradient of $V(\textbf{x})$; for instance, the partial derivative of $V(\textbf{x})$ with respect to $x_E$ is as follows 
\begin{align}
 \left.
	 \begin{array}{l l}
	\frac{\partial V(\textbf{x})}{\partial x_E} \!\!\!\!&= \frac{\nu_E x_E - \frac{1}{2}[ (y_2-y_3) R_1^2 + (y_3-y_1) R_2^2  + (y_1-y_2) R_3^2]  }{\Lambda}  \\
	&~~ -  \mu_E \frac{ \nu_E R_E^2 - \nu_1 R_1^2  -\nu_2 R_2^2  -\nu_3 R_3^2  }{2\Lambda^2}.
\end{array}  \right.  \nonumber 
\end{align}
Adding and subtracting the term $\nu_E \frac{ x^*}{\Lambda}$  to the previous equation
and rearranging terms we obtain the following 
\begin{align}
 \left.
	 \begin{array}{l l}
	\frac{\partial V(\textbf{x})}{\partial x_E} \!\!\!\!&= \frac{\nu_E (x_E-x^*) - \frac{1}{2}[ (y_2-y_3) R_1^2 + (y_3-y_1) R_2^2  + (y_1-y_2) R_3^2]  }{\Lambda}  \\
	&~~ +\nu_E \frac{x^*}{\Lambda} -  \mu_E \frac{ \nu_E R_E^2 - \nu_1 R_1^2  -\nu_2 R_2^2  -\nu_3 R_3^2  }{2\Lambda^2} \\
	& = \frac{\nu_E (x_E-x^*) - \frac{1}{2}[ (y_2-y_3) R_1^2 + (y_3-y_1) R_2^2  + (y_1-y_2) R_3^2]  }{\Lambda}  \\
	&~~ + \frac{ ( \mu_E \nu_1- \nu_E \mu_1 )R_1^2  +(\mu_E\nu_2 - \nu_E\mu_2)R_2^2  +(\mu_E\nu_3 - \nu_E\mu_3) R_3^2  }{2\Lambda^2}
\end{array}  \right.  \nonumber 
\end{align}
where $x^* =\frac{\mu_E R_E^2 - \mu_1 R_1^2  -\mu_2 R_2^2  -\mu_3 R_3^2}{2\Lambda}$, as given by \eqref{eq:Xopt}, was substituted in the second line of the previous equation. Finally, grouping like terms together we have that
\begin{align}
 \left.
	 \begin{array}{l l}
	\frac{\partial V(\textbf{x})}{\partial x_E} \!\!\!\!&= \frac{\nu_E (x_E-x^*)}{\Lambda} + \frac{  \mu_E \nu_1- \nu_E \mu_1 - (y_2-y_3) \Lambda }{2\Lambda^2} R_1^2 \\
	&~~+ \frac{ \mu_E\nu_2 - \nu_E\mu_2 - (y_3-y_1) \Lambda }{2\Lambda^2} R_2^2  \\
	&~~+ \frac{ \mu_E\nu_3 - \nu_E\mu_3 - (y_3-y_1) \Lambda }{2\Lambda^2} R_3^2. 
\end{array}  \right.  \nonumber 
\end{align}
It can be shown that 
\begin{align}
 \left.
	 \begin{array}{l l}
	  \mu_E \nu_1- \nu_E \mu_1 - (y_2-y_3) \Lambda = 0 \\
	 \mu_E\nu_2 - \nu_E\mu_2 - (y_3-y_1) \Lambda =0  \\
	 \mu_E\nu_3 - \nu_E\mu_3 - (y_3-y_1) \Lambda = 0. 
\end{array}  \right.  \nonumber 
\end{align}
Hence, the partial derivative of $V(\textbf{x})$ with respect to $x_E$ is given by
\begin{align}
 \left.
	 \begin{array}{l l}
	\frac{\partial V(\textbf{x})}{\partial x_E} \!\!\!\!&= \frac{\nu_E }{\Lambda}  (x_E-x^*).
\end{array}  \right.  \nonumber 
\end{align}
The partial derivatives $\frac{\partial V(\textbf{x})}{\partial y_E}$ and $\frac{\partial V(\textbf{x})}{\partial z_E}$ can be obtained following similar steps and we have that the partial derivative of $V(\textbf{x})$ with respect to the state of the evader $\textbf{x}_E$ is as follows
\begin{align}
 \left.
	 \begin{array}{l l}
	\frac{\partial V(\textbf{x})}{\partial \textbf{x}_E} \!\!\!\!&= \frac{\nu_E }{\Lambda} [x_E-x^*(\textbf{x}),  \ y_E - y^*(\textbf{x}), \ z_E-  z^*(\textbf{x}) ].
\end{array}  \right.  \label{eq:PVPxE}
\end{align}
Likewise, the partial derivatives of $V(\textbf{x})$ with respect to the state of each pursuer are given by
\begin{align}
 \left.
	 \begin{array}{l l}
	\frac{\partial V(\textbf{x})}{\partial \textbf{x}_i} \!\!\!\!&= - \frac{\nu_i }{\Lambda} [x_i-x^*(\textbf{x}),  \ y_i - y^*(\textbf{x}), \ z_i-  z^*(\textbf{x}) ]  \\
\end{array}  \right.  \label{eq:PVPx1}
\end{align}
for $i=1,2,3$. Note that, for $\textbf{x} \in \mathcal{R}_p$, we have that $\Lambda \neq 0$, that is, the three planes given by \eqref{eq:OP-Ei} intersect at only one point. Therefore, the Value function $V(\textbf{x})$ is continuous and continuously differentiable for any $\textbf{x} \in \mathcal{R}_p$.

Finally, we show that $V(\textbf{x})$ is the solution of the HJI equation.
The HJI equation  is given in general by
$-\frac{\partial V(t,\textbf{x})}{\partial t} =\frac{\partial V(t,\textbf{x})}{\partial \textbf{x}}\cdot  \textbf{f}(\textbf{x},\textbf{u}_E^*,\textbf{u}_i^*) + g(t,\textbf{x},\textbf{u}_E^*,\textbf{u}_i^*) $.
In this problem we have $\frac{\partial V(t,\textbf{x})}{\partial t}=0$,  $g(t,\textbf{x},\textbf{u}_E^*,\textbf{u}_i^*)=0$ and we have that
\begin{align}
 \left.
	 \begin{array}{l l}
	  \frac{\partial V(\textbf{x})}{\partial \textbf{x}} \! \cdot \! \textbf{f}(\textbf{x},\textbf{u}_E^*,\textbf{u}_i^*) 
	  =  - \frac{\nu_E[(x^*-x_E)^2+(y^*-y_E)^2+(z^*-z_E)^2 ] }{\Lambda d_E}  \\  
 \qquad \qquad \qquad \qquad	\ \  + \frac{\nu_1[(x^*-x_1)^2+(y^*-y_1)^2+(z^*-z_1)^2 ] }{\Lambda d_1}  \\  
  \qquad \qquad \qquad \qquad	\ \  + \frac{\nu_2[(x^*-x_2)^2+(y^*-y_2)^2+(z^*-z_2)^2 ] }{\Lambda d_2}  \\  
   \qquad \qquad \qquad \qquad	 \ \  + \frac{\nu_3[(x^*-x_3)^2+(y^*-y_3)^2+(z^*-z_3)^2 ] }{\Lambda d_3}  \\
    \qquad \qquad \qquad \qquad = - \frac{\nu_E d_E }{\Lambda} +  \frac{\nu_1 d_1 }{\Lambda} + \frac{\nu_2 d_2 }{\Lambda} + \frac{\nu_3 d_3 }{\Lambda}.
\end{array}  \right.   \nonumber
\end{align}	
Note that the point $I^*$ is equidistant to the locations of all players, then $d_E=d_i$, thus, we can write
\begin{align}
 \left.
	 \begin{array}{l l}
	  \frac{\partial V(\textbf{x})}{\partial \textbf{x}}\cdot  \textbf{f}(\textbf{x},\textbf{u}_E^*,\textbf{u}_i^**) = 
	   \frac{d_E }{\Lambda} [ \nu_1 +\nu_2 +\nu_3 -\nu_E] .
\end{array}  \right.   \nonumber
\end{align}
Using \eqref{eq:zoterms}, it is easy to verify that $\nu_1 +\nu_2 +\nu_3 -\nu_E=0$. Hence, $\frac{\partial V(\textbf{x})}{\partial \textbf{x}}\cdot  \textbf{f}(\textbf{x},\textbf{u}_E^*,\textbf{u}_i^**) = 0$.
In conclusion, we have shown that the $C^1$ Value function \eqref{eq:ValueFn} is the solution of the HJI equation and the associated state-feedback strategies \eqref{eq:OptimalInputE}-\eqref{eq:OptimalInputsP} are the optimal strategies of the differential game. $\square$

\textit{Remark}. The results in this section can be extended to consider groups of more than three pursuers. In such a case the reachable region of the evader is determined by the intersection of several half-spaces and the evader's optimal aimpoint is given by the reachable vertex which is the closest to $\Omega_G$. In general, only three pursuers are active, that is, only three pursuers (those three pursuers associated with the optimal aimpoint which is intersection point of three planes) will eventually capture the evader under optimal play. Also note that in the presence of only one or only two pursuers, the evader is always able to reach $\Omega_G$ and win the game. The configurations analyzed in \cite{yan2019construction} where one or two pursuers seem to be the winners are actually singular surfaces of dispersal type under the more general class of state-feedback strategies. Following the classical work by Isaacs \cite{Isaacs65} one can easily see that the pursuers are never able to correctly guess the strategy of the evader and the plane (in the case of one pursuer) or the intersecting line (in the case of two pursuers) will inevitably tilt and provide a route to the evader in order to reach $\Omega_G$. Therefore, at least three pursuers are needed in order to guard a plane in the 3-D space.

\begin{figure}
	\begin{center}
		\includegraphics[width=8.4cm,trim=1.4cm .6cm 2.2cm .7cm]{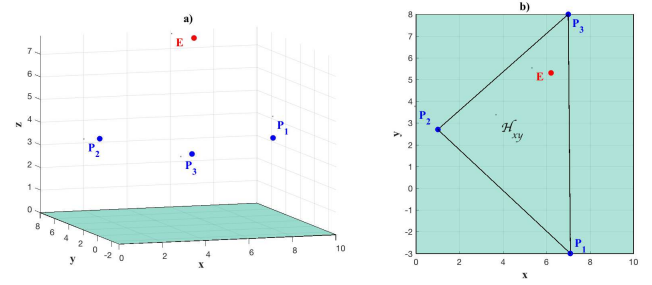}
	\caption{Game of guarding a plane. a) Evader and three pursuers in the 3-D space with the plane $z(x,y)=0$ shown in color green. b) Players seen from above with the $x$-$y$ plane projection of the pursuers' convex hull}
	\label{fig:ConvexHull}
	\end{center}
\end{figure}

\section{The Barrier Surface} \label{sec:Barrier}
We now solve the Game of Kind in closed-form in order to characterize the Barrier surface which separates the winning regions of the players. The Barrier surface is obtained when the evader is intercepted by the pursuers at exactly the same time instant when he reaches the goal plane. Let us define 
\begin{align}
 \left.
	 \begin{array}{l l}
	a(\textbf{x}) = \frac{1}{2} \frac{(y_2-y_3)R_1^2 + (y_3-y_1)R_2^2 + (y_1-y_2)R_3^2}{x_1(y_2-y_3)+x_2(y_3-y_1)+x_3(y_1-y_2)} \\
	b(\textbf{x}) = \frac{1}{2} \frac{(x_3-x_2)R_1^2 + (x_1-x_3)R_2^2 + (x_2-x_1)R_3^2}{x_1(y_2-y_3)+x_2(y_3-y_1)+x_3(y_1-y_2)}.
\end{array}  \right.   \label{eq:aPar}
\end{align}
Let us also consider the projected convex hull of the pursuers onto the $x$-$y$ plane. Fig. \ref{fig:ConvexHull}.b provides an illustration of the projected convex hull, that is, as it is seen from `above'. In other words, the $z$-coordinate of the players is disregarded. Let $\mathcal{H}_{xy}$ denote the projected convex hull of the pursuers onto the $x$-$y$ plane. Define
\begin{align}
 \left.
	 \begin{array}{l l}
	\mathcal{B}_{xy}=  \big\{ \ \textbf{x} \ | \ z_E>0, \ (x_E,y_E) \in \mathcal{H}_{xy}  \big\}. 
\end{array}  \right.   \nonumber
\end{align}
The following theorem provides the solution of the Game of Kind.
\begin{theorem}
Consider the differential game of protecting a plane in the 3-D space. The Barrier surface is given by $\mathcal{B}=\mathcal{B}_{xy}  \cap \mathcal{B}_s$ where 
\begin{align}
 \left.
	 \begin{array}{l l}
	\mathcal{B}_s =  \big\{ \ \textbf{x} \ | \ B( \textbf{x})=0  \big\} 
\end{array}  \right.   \nonumber
\end{align}
and the Barrier function is given by
\begin{align}
 \left.
	 \begin{array}{l l}
	B( \textbf{x})\!\!\!& = (x_E-a(\textbf{x}))^2 + (y_E-b(\textbf{x}))^2 +z_E^2  \\
	&~~-  (x_1-a(\textbf{x}))^2 - (y_1-b(\textbf{x}))^2 - z_1^2.
\end{array}  \right.   \label{eq:BarrFn}
\end{align}
\end{theorem}
\textit{Proof}.
The Barrier surface is obtained when $z^*(\textbf{x})=0$. Since $z^*(\textbf{x})=V(\textbf{x})$, then we set \eqref{eq:ValueFn} equal to zero as follows
\begin{align}
 \left.
	 \begin{array}{l l}
	 \frac{1}{\Lambda}\big( [x_1(y_2\!-\!y_3)+x_2(y_3\!-\!y_1)+x_3(y_1\!-\!y_2)]R_E^2  \\
	 ~~  -[x_2(y_3\!-\!y_E)+x_3(y_E\!-\!y_2)+x_E(y_2\!-\!y_3)]R_1^2  \\ 
	 ~~  -[x_3(y_1\!-\!y_E)+x_1(y_E\!-\!y_3)+x_E(y_3\!-\!y_1)]R_2^2  \\
	 ~~  -[x_1(y_2\!-\!y_E)+x_2(y_E\!-\!y_1)+x_E(y_1\!-\!y_2)]R_3^2   \big) =0.
\end{array}  \right.  \nonumber
\end{align}
We now multiply both sides of the previous equation by $\Lambda$. Also, we divide both sides of the equation by the term $x_1(y_2\!-\!y_3)+x_2(y_3\!-\!y_1)+x_3(y_1\!-\!y_2)$ to obtain
\begin{align}
 \left.
	 \begin{array}{l l}
	x_E^2+y_E^2+z_E^2  - \frac{x_2(y_3-y_E)+x_3(y_E-y_2)+x_E(y_2-y_3)}{x_1(y_2-y_3)+x_2(y_3-y_1)+x_3(y_1-y_2)} R_1^2  \\ 
	  -\frac{x_3(y_1-y_E)+x_1(y_E-y_3)+x_E(y_3-y_1)}{x_1(y_2-y_3)+x_2(y_3-y_1)+x_3(y_1-y_2)} R_2^2  \\
	  -\frac{x_1(y_2-y_E)+x_2(y_E-y_1)+x_E(y_1-y_2)}{x_1(y_2-y_3)+x_2(y_3-y_1)+x_3(y_1-y_2)}  R_3^2    =0.
\end{array}  \right.  \nonumber
\end{align}
Rearranging terms in the previous equation we have that
\begin{align}
 \left.
	 \begin{array}{l l}
	~~ x_E^2 - \frac{(y_2-y_3)R_1^2 + (y_3-y_1)R_2^2 + (y_1-y_2)R_3^2}{x_1(y_2-y_3)+x_2(y_3-y_1)+x_3(y_1-y_2)} x_E  \\ 
	+y_E^2  - \frac{(x_3-x_2)R_1^2 + (x_1-x_3)R_2^2 + (x_2-x_1)R_3^2}{x_1(y_2-y_3)+x_2(y_3-y_1)+x_3(y_1-y_2)} y_E\\
	+z_E^2  - \frac{(x_2y_3 -x_3y_2)R_1^2 + (x_3y_1-x_1y_3)R_2^2 + (x_1y_2 -x_2y_1)R_3^2}{x_1(y_2- y_3)+x_2(y_3 -y_1)+x_3(y_1-y_2)} =0 \\
	\Rightarrow ~~ x_E^2 -2x_E a(\textbf{x}) + y_E^2 -2y_Eb(\textbf{x}) + z_E^2  \\
	 \qquad - \frac{(x_2y_3 -x_3y_2)R_1^2 + (x_3y_1-x_1y_3)R_2^2 + (x_1y_2 -x_2y_1)R_3^2}{x_1(y_2- y_3)+x_2(y_3 -y_1)+x_3(y_1-y_2)} =0. 
	  \end{array}  \right.  \nonumber
\end{align}
Adding and subtracting the terms $a^2(\textbf{x})$ and $b^2(\textbf{x})$ we obtain the following
\begin{align}
 \left.
	 \begin{array}{l l}
	(x_E -a(\textbf{x}))^2 + (y_E -b(\textbf{x}))^2 + z_E^2 - a^2(\textbf{x}) -b^2(\textbf{x})\\
	  - \frac{(x_2y_3 -x_3y_2)R_1^2 + (x_3y_1-x_1y_3)R_2^2 + (x_1y_2 -x_2y_1)R_3^2}{x_1(y_2- y_3)+x_2(y_3 -y_1)+x_3(y_1-y_2)} =0 \\
	  \Rightarrow ~~	(x_E -a(\textbf{x}))^2 + (y_E -b(\textbf{x}))^2 + z_E^2 \\
	\qquad -(x_1-a(\textbf{x}))^2 - (y_1-b(\textbf{x}))^2 - z_1^2 = 0.
	  \end{array}  \right.  \nonumber
\end{align}
Hence, the Barrier surface is obtained when $B( \textbf{x})=0$, where $B( \textbf{x})$ is given by \eqref{eq:BarrFn}. \ $\square$

The Barrier surface is such that $\mathcal{B}\in\mathbb{R}^{12}$, that is, $\mathcal{B}$ is a surface within the game's state space of dimension twelve. However, an informative cross-section of the Barrier surface can be obtained by fixing the position of the pursuers. Therefore, it is possible to plot the the cross-section of the Barrier surface in the 3-D space, in terms of the potential position of the evader. This process provides a clear illustration of the players' winning regions in the 3-D space. 

\begin{corollary}
Let $\mathcal{C}\in\mathbb{R}^3$ denote the cross-section of the Barrier surface  when the pursuers' positions are fixed.
The cross-section of the Barrier surface is given by $\mathcal{C}= \mathcal{C}_{xy} \cap \mathcal{C}_s$, where
\begin{align}
 \left.
	 \begin{array}{l l}
	\mathcal{C}_{xy} =\{ x,y,z  \ | \ z>0, \  (x,y) \in \mathcal{H}_{xy} \}, 
\end{array}  \right.   \nonumber
\end{align}	
\begin{align}
 \left.
	 \begin{array}{l l}	
	\mathcal{C}_s =\{ x,y,z  \ | &(x-a)^2 + (y-b)^2  +z^2  \\
	& =  (x_1-a)^2 + (y_1-b)^2 + z_1^2 \},
\end{array}  \right.   \nonumber
\end{align}
and the parameters $a$ and $b$ are given by \eqref{eq:aPar}, with fixed pursuers' positions. \ $\square$
\end{corollary}

It is easy to see that $\mathcal{C}_s$ is the sphere with center at $(a,b,0)$ and radius $r=\sqrt{(x_1-a)^2 + (y_1-b)^2 + z_1^2}$. This expression provides a simple and explicit form of the Barrier surface.
Compared to \cite{yan2019construction}, where only preliminary conditions on the state are used for the Barrier surface, in this paper we have obtained a closed-form solution of the Barrier surface, that is, the Barrier surface is completely characterized as a sphere where the center coordinates and the radius are given in simple and explicit form.

\subsection{Singular Surface and Future Work}
 A challenging situation for a group of pursuers trying to protect a plane of infinite length from a same speed evader is that the evader may try to escape the projected convex hull. Let $\partial \mathcal{H}_{xy}$ be the boundary of $\mathcal{H}_{xy}$. When $\textbf{x} \in \partial \mathcal{H}_{xy}$ the pursuers should not let the evader escape from the projected convex hull. This strategy requires a guidance parallel to the plane by the pursuers to prevent such outcome. If all players keep such strategy, then a non-termination outcome will occur and the distance between the evader and the goal plane remains constant. However, by Definition 1, non-termination is an inferior outcome for the evader and the pursuers will benefit from extending the duration of the game. 
 
Designing strategies on singular surfaces is a complex and difficult process which generally depends on the particular game on hand. The pursuers in this problem will need to switch their guidance between aiming at the intersection point of the three half-planes and implementing the parallel strategy to keep the evader within the projected convex hull. Since the evader is free to implement non-optimal strategies, frequent switching by the pursuers may occur. In order to avoid frequent pursuers' switching behaviors, mixed and/or degenerate strategies may provide an appropriate solution which also minimize risk with respect to uncertain behaviors by the evader. 
 It is important to highlight that all these strategies are state-feedback strategies, that is, the players never share strategic information with the opponent. These behaviors will be exemplified in the following section.

There are several extensions to this problem, which are interesting on their own, but they also provide significantly different outcomes and types of Barrier surface. 

\textit{Positive capture radius}. A future extension will address the case where the pursuers are endowed with a positive capture radius, that is, point capture is not necessary and the evader is intercepted if it is within distance of a given pursuer $i$ equal to that pursuer's capture radius. In such a case, the evader cannot extend indefinitely the duration of the game without increasing its separation with respect to the goal plane.

\textit{Players with different speeds}. In the case where the pursuers are faster than the evader, the evader cannot avoid termination of the game regardless of whether point capture or positive capture distance is assumed. The analysis of reachable regions between the evader and any given pursuers is carried over by using Apollonius spheres instead of planes.

On the other hand, the case where the evader is faster than the pursuers is a much more complex scenario where different types of strategies can be performed by the evader in order to circumvent the pursuers and reach the goal set. In such a case, the pursuers need to be endowed with a positive capture radius, otherwise it is no possible for a slow pursuer to point capture a fast evader \cite{Garcia21auto}.

\textit{Guarding a finite size plane with same speed players and point capture}. Finally, a very practical situation is concerned with the protection of a plane of finite dimensions. This problem represents a real-world application where, for instance, a team of players is tasked to protect a finite length coastline and its airspace from intruders. Besides interception of $E$, termination of the game is also achieved when $E$ is chased out of the game set, \textit{i.e.}, the boundary of the coastline. In this case, point capture and same speed players can be considered and the previously mentioned singular surface is not present. In a similar fashion, optimal strategies for protection of other target sets, such as a line or a disk on a plane can be obtained using similar ideas and concepts and will be addressed in the future along with their necessary conditions, \textit{i.e.} the number of pursuers necessary to protect these type of targets. 

A common feature among all three scenarios is that the shape of the Barrier surface is significantly different from the Barrier surface obtained in this paper, since the winning region of the pursuers is considerably expanded (each in its own way) beyond the projected convex hull of the pursuers. In other words, switching strategies are not necessary by the pursuers since the evader can be located outside the projected convex hull and still be unable to win the game.
It is interesting and useful to exactly characterize the Barrier surface in each case in order to accurately delineate the winning regions of each player or team. These problems will be addressed in future work.

\begin{figure}
	\begin{center}
		\includegraphics[width=8.2cm,trim=1.5cm .7cm 1.4cm .6cm]{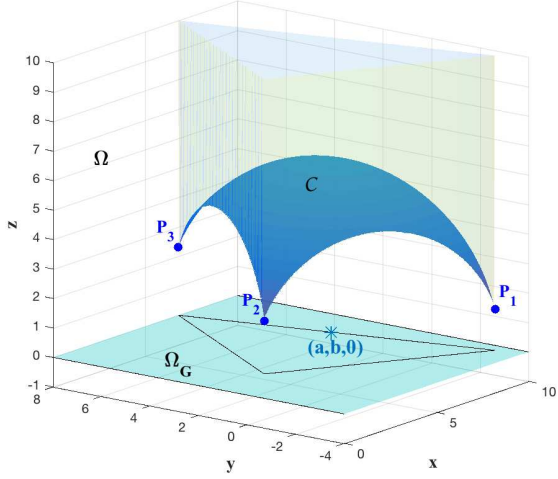}
	\caption{Example of the cross-section of the Barrier surface showing the projection of the pursuers' convex hull and the location of the center of the sphere $(a,b,0)$}
	\label{fig:Sphere}
	\end{center}
\end{figure}

\begin{figure}
	\begin{center}
		\includegraphics[width=8.2cm,trim=1.5cm 1.1cm 1.4cm 1.4cm]{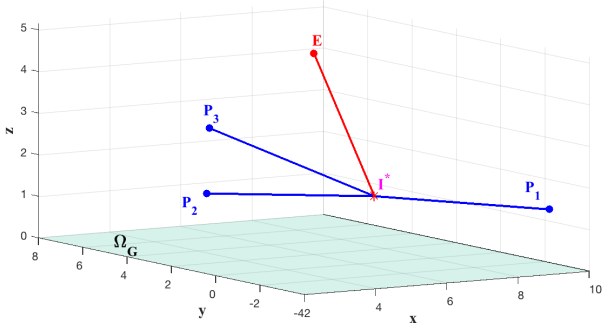}
	\caption{Optimal trajectories in Example 2 where the evader is intercepted by the group of pursuers at the optimal interception point $I^*=(x^*,y^*,z^*)$}
	\label{fig:Traj}
	\end{center}
\end{figure}

\begin{figure}
	\begin{center}
		\includegraphics[width=7.2cm,trim=1.5cm -.5cm 1.4cm .1cm]{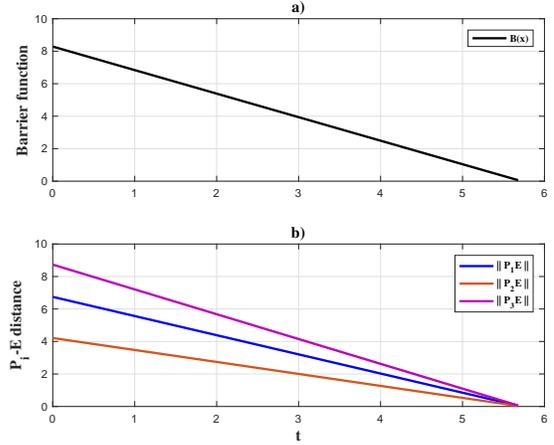}
	\caption{Example 2: a) Barrier function. b) Distance between evader and each pursuer}
	\label{fig:ex22}
	\end{center}
\end{figure}

\section{Examples} \label{sec:ex}
\textit{Example 1. Barrier surface}. Consider the initial positions of the pursuers as follows, $P_1=(9.5,  -3, 1.4)$, $P_2=(2.5, 1.2, 1.8)$,  $P_3=(6.8, 8, 2.3)$. The cross-section of the Barrier surface, $\mathcal{C}$, is shown in Fig. \ref{fig:Sphere}. It is a section of the sphere, delineated by the projection of the pursuers' convex hull, with center at $(a,b,0)$, where $a=8.020$ and $b=2.619$, and radius $r=\sqrt{(x_1-a)^2 + (y_1-b)^2 + z_1^2}=5.977$.

\textit{Example 2. Optimal strategies}. Consider the same initial positions of the pursuers as in Example 1. In addition, let the initial position of the evader be given by $E=(4.76, 0, 5.15)$ and we have that $B(\textbf{x}_0)>0$. Therefore, the pursuers are  able to win the game if they follow their optimal trajectories. Fig. \ref{fig:Traj} shows the optimal trajectories when all players, the evader and all pursuers, implement their optimal strategies \eqref{eq:OptimalInputE}-\eqref{eq:OptimalInputsP}. The evader is synchronously captured by all pursuers and the optimal interception point is given by $I^*=(8.032, 2.524,  1.195)$. The Value of the game is $V(\textbf{x})=z_E(t_f)=1.195$. In Fig. \ref{fig:ex22}.a the Barrier function, $B(\textbf{x})$ is shown as a function of time. It can be seen that $B(\textbf{x})>0$ for $t\in[0,t_f)$. Fig.  \ref{fig:ex22}.b shows the separation between the evader and each one of the pursuers. It can be seen that all separations become equal to zero at the terminal time $t_f$.

\begin{figure}
	\begin{center}
		\includegraphics[width=8.2cm,trim=1.5cm 1.1cm 1.5cm .8cm]{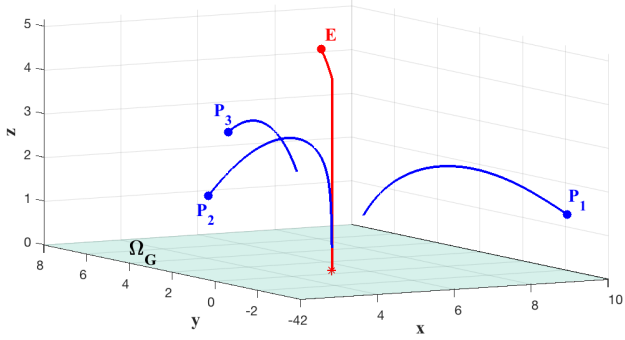}
	\caption{Optimal trajectories in Example 3 where the pursuers play non-optimally and implement the pure pursuit strategy}
	\label{fig:Ex3Traj}
	\end{center}
\end{figure}

\begin{figure}
	\begin{center}
		\includegraphics[width=7.8cm,trim=1.5cm .0cm 1.4cm .0cm]{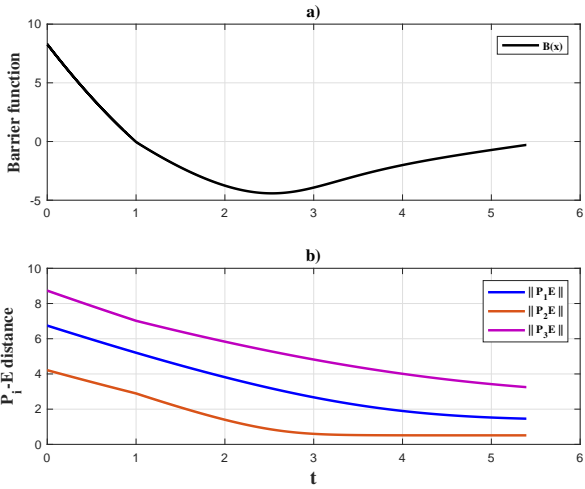}
	\caption{Example 3: a) Barrier function. b) Distance between evader and each pursuer}
	\label{fig:ex32}
	\end{center}
\end{figure}

\textit{Example 3. Non-optimal play by pursuers}. Consider the same initial conditions as in Example 2. We now compare the outcome and performance to the case where the pursuers do not implement their optimal strategy and use Pure Pursuit (PP) guidance on the evader. The PP guidance is a common strategy and it is widely used in many pursuit-evasion scenarios. However, PP is not the optimal strategy for this differential game. The resulting trajectories are shown in Fig. \ref{fig:Ex3Traj}, where it can be seen that the evader is actually able to win the game by reaching $\Omega_G$ before being intercepted by any of the pursuers. $B(\textbf{x})$ as a function of time is shown in Fig. \ref{fig:ex32}.a. $B(\textbf{x})$ changes sign at time $t=0.994$. Since $B(\textbf{x})$ can be easily computed from \eqref{eq:BarrFn}, the evader, who is implementing its optimal strategy, is able to check the sign of $B(\textbf{x})$. The evader implements the state-feedback optimal strategy \eqref{eq:OptimalInputE} in closed-loop manner and then it switches guidance and it heads directly to $\Omega_G$ once $B(\textbf{x})<0$. Fig. \ref{fig:ex32}.b shows the separations between the evader and each one of the pursuers. None of the pursuers is able to capture $E$. The pursuers perform very poorly by disregarding their optimal strategy and choosing a non-optimal guidance such as PP.

\begin{figure}
	\begin{center}
		\includegraphics[width=8.2cm,trim=1.5cm 2.2cm 1.5cm 1.2cm]{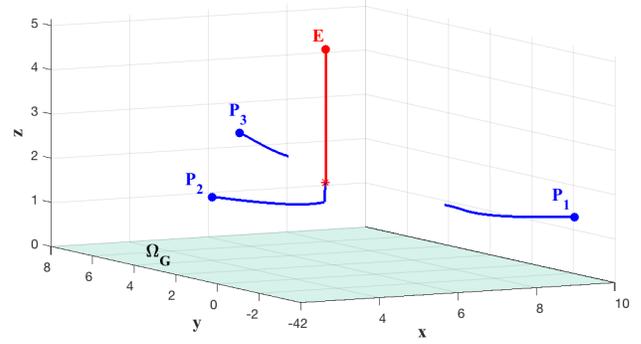}
	\caption{Optimal trajectories in Example 4 where the evader plays non-optimally}
	\label{fig:Ex4Traj}
	\end{center}
\end{figure}

\begin{figure}
	\begin{center}
		\includegraphics[width=7.8cm,trim=1.5cm .1cm 1.4cm .0cm]{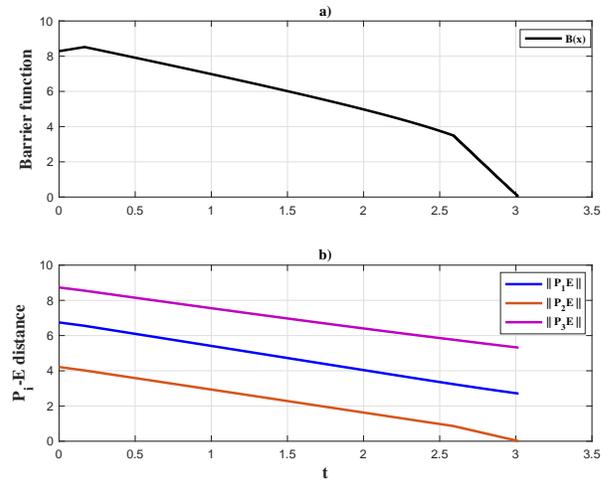}
	\caption{Example 4: a) Barrier function. b) Distance between evader and each pursuer}
	\label{fig:ex42}
	\end{center}
\end{figure}

\textit{Example 4. Non-optimal play by evader}. Finally, we consider the opposite case to Example 3; now, the evader plays non-optimally while the pursuers implement the state-feedback optimal strategy \eqref{eq:OptimalInputsP} in closed-loop manner. The evader, in this case, implements a greedy strategy and tries to reach $\Omega_G$ by heading directly into the goal plane and disregarding its optimal strategy and the Game of Kind solution. However, this strategy is not known nor assumed by the pursuers who need to implement strategies based only on state measurements. 
In this example, the pursuers initially aim at the optimal interception point. As the pursuers move, the convex hull changes and they realize that $E$'s position is close to $\partial \mathcal{H}_{xy}$, along the $P_1-P_2$ segment. The corresponding pursuers, $P_1$ and $P_2$, implement a mixed strategy and aim at the intersection of the three planes $H_1$, $H_2$, and the plane orthogonal to the $z$-plane along the $P_1$ and $P_2$ positions (or the wall along the $P_1-P_2$ segment of the convex hull). In this way they protect the convex hull from the external side while $P_3$ closes in from the internal side.
A further switch to parallel guidance will be needed in the case $E$ tries to escape the convex hull and a non-termination outcome will follow. Such is not the case in this example since $E$ aims directly at $\Omega_G$. As the pursuers keep moving and $E$ keeps implementing a non-optimal strategy, a second guidance switch is implemented by the pursuers since now $E$ is close to $\partial \mathcal{H}_{xy}$ but now along the $P_2$ vertex. $P_2$ now aims directly at $E$ while the remaining pursuers protect the convex hull from either side each. The evader is intercepted by only one pursuer, as it is typical when $E$ does not implement its optimal guidance.

The terminal separation is $z_E(t_f)=2.130> V(\textbf{x})$. The resulting trajectories are shown in Fig. \ref{fig:Ex4Traj}. Not only $E$ was captured, but the terminal separation was significantly increased with resect to the Value of the game. This is good for the pursuers and it is a consequence of $E$ not playing optimally. Figure \ref{fig:ex42}.a shows the Barrier function. The switches on the pursuers' guidance occur at $t=0.169$ and at $t=2.588$ and they can be noticed on the sharp corners of $B(\textbf{x})$. Figure  \ref{fig:ex42}.b shows the separation between the evader and each one of the pursuers. It can be seen that the separation between $P_2$ and $E$ becomes equal to zero at the terminal time $t_f$. We emphasize that all the strategies implemented by the pursuers are state-feedback strategies. They employ a high level of cooperation in order to intercept $E$ while blocking escape routes under uncertainty of $E$'s strategy.

\section{Conclusions} \label{sec:concl}
The game of protecting a plane in the 3-D space was addressed in this paper. Significant extensions with respect to recent work were presented. First, the more general and practical class of state-feedback strategies was utilized. Second, the Barrier surface was characterized in a simple and explicit form. Finally and more importantly, the solution of the Game of Degree was obtained, that is, the optimal strategies of each player were obtained and verified. The last point represents a valuable contribution since, only by solving the Game of Degree, the players are actually able to achieve the prescribed outcome of the Game of Kind.

\bibliographystyle{IEEEtran}
\bibliography{ReferencesTAD}

\end{document}